\documentclass[11pt]{paper}

\usepackage{amssymb,amsmath,enumerate,theorem,epsfig,rotating,graphics,changebar,eepic}\usepackage{amsfonts}
\usepackage{a4,latexsym,parskip}
\usepackage{hyperref}
\hypersetup{pdfauthor=MSAS} \hypersetup{pdftitle=DS}

\newcommand{\PP}{\mathbb{P}}

\newcommand{\C} [1][]{\mathbb{C}^{#1}}
\newcommand{\Q} [1] []{\mathbb{Q}_{#1}}
\newcommand{\N} [1][] {\mathbb{N}_{#1}}

\newcommand{\F}{\mathbb{F}}
\newcommand{\Z}{\mathbb{Z}}

\newcommand{\NN}{\mathcal{N}}

\newcommand{\qed}{\hfill \ensuremath{\Box}}

\newtheorem{Spezial-Theorem}{Theorem}[section]
\newtheorem{Spezial-Proposition}{Proposition}[section]
\theoremstyle{break} \newtheorem{Theorem}{Theorem}[section]
\newtheorem{Proposition}[Theorem]{Proposition}
\newtheorem{Lemma}[Theorem]{Lemma}
\newtheorem{Example}[Theorem]{Example}

\newtheorem{Corollary}[Theorem]{Corollary}

\newtheorem{Remark}[Theorem]{Remark}

\newtheorem{Criterion}[Theorem]{Criterion}

\begin{document}
\setlength{\unitlength}{1cm}

\title{On Davenport-Stothers inequalities and elliptic surfaces in positive characteristic}
\author{Matthias Sch\"utt, Andreas Schweizer}

%\date{\today}

\maketitle

\abstract{We show that the Davenport-Stothers inequality from characteristic
$0$ fails in any characteristic $p>3$. The proof uses elliptic surfaces over
$\PP^1$ and inseparable base change. We then present adjusted inequalities. These follow from results of Pesenti-Szpiro. For characteristic $2$ and $3$, we achieve a similar
result in terms of the maximal singular fibres of elliptic surfaces
over $\PP^1$. Our ideas are also related to supersingular
surfaces (in Shioda's sense).}



\vspace{0.2cm}

\section{Introduction}
\label{s:intro}

We recall a result of Davenport:

\begin{Theorem}[Davenport {\cite{Dav}}]\label{Thm:DS}
Let $M\in\N$  and $f,g \in\mathbb{C}[t]$ with deg
$f=2M$, deg $g=3M$. Then
\[ \text{deg}(f^3-g^2)\geq M+1 \Leftrightarrow f^3\neq g^2.
\]
\end{Theorem}
Later Stothers \cite{St} showed among other results that for every
$M$ there do indeed exist polynomials $f$ and $g$ with
$\text{deg}(f^3 -g^2)=M+1$, and that modulo affine transformations of
$t$ and making $f$ and $g$ monic their number is finite for each $M$.
\par
An extensive account of this and interesting applications
to elliptic surfaces can be found in Shioda's paper \cite{Sh-DS}.
Shioda remarks that the inequality holds in characteristic $p$
(i.e.~for polynomials $f,g\in k[t]$ over an algebraically closed field $k$ of characteristic $p$) if $p>6M$. We will
refer to this as \emph{$DS(M)$ holds mod $p$}. For $M=3,4$,
Shioda's results were complemented in \cite{S-max} to prove that
$DS(M)$ holds mod $p$ for any $p>3$. (For $M=1,2$, this follows from the theory of rational elliptic surfaces.)

The biggest part of this paper is an extensive study of $DS(M)$ mod $p$ 
for general $M$. The main results are collected in the following

\begin{Theorem}\label{Thm:p}
Let $k$ be an algebraically closed field of characteristic $p>3$. 
Let $f, g\in k[t]$ with $\deg f=2M$, $\deg g=3M$ for some $M\in\N$.

\begin{enumerate}[(a)]
\item For every $p>3$ there are infinitely many $f$ and $g$ in $\F_p[t]$
with $f^3 -g^2 =-(t^2 -1)^2$. In particular, there are such examples with
arbitrarily large $M$.

\item Let $(f,g)$ be a counterexample to $DS(M)$ mod $p$. Then every factor 
of $f^3 -g^2$ has multiplicity at least $2$. If such a factor does not 
divide $f$ and $g$, its multiplicity is divisible by $p$. 

\item In every counterexample to $DS(M)$ mod $p$, the polynomials $f$ and 
$g$ have at least two different common zeroes; in other words,
\[ 
\# \{s\in k; f(s)=g(s)=0\}\leq 1 \Rightarrow \text{deg}(f^3-g^2)\geq M+1.
\]
In particular, the Davenport-Stothers inequality holds for coprime
$f$ and $g$.

\item Let $M>2$. Then the following relation holds:
\[ 
f^3\neq g^2 \Leftrightarrow \text{deg}(f^3-g^2)\geq \begin{cases} 4 & \text{if $M$ is odd,}\\
5 & \text{if $M$ is even, } p\equiv 7\text{ mod } 12\;\; \text{ or } M=4,\\
6 & \text{if $M$ is even, } M>4, \;p\not\equiv 7\text{ mod }12.
\end{cases}
\]
In general, these bounds are sharp.

\item $DS(M)$ holds mod $p$ if there is no $n\in\N$ such that
\[
5M\leq np\leq 6M-\begin{cases} 4, & \text{if $M$ is odd,}\\
5, & \text{if $M$ is even.} \end{cases}
\]
In particular, this holds true if $p>6M-5$.

\item  For every $p>3$ there is an $M_0\in\N$ such that $DS(M)$ does not 
hold mod $p$ for any $M\geq M_0$. For $p>29$, the optimal value is
\[
M_0 = \begin{cases} \frac{5p+7}6, & \text{if} \;\; p\equiv 1\mod 6,\\
p+2, & \text{if}\;\; p\equiv 5\mod 6.\end{cases}
\]

\item For any $M>4$, there is a $p>3$ such that $DS(M)$ does not hold mod $p$.
Moreover, a counterexample can be given over $\F_p$.
\end{enumerate}

\end{Theorem}

These results will be proved, partially with slightly more detail information,
as Example \ref{Ex:deg4}, Remark \ref{Rem:multiplicity}, 
Proposition \ref{Lem:coprime}, Proposition \ref{Prop:4,5}, 
Criterion \ref{Crit}, Corollary \ref{Cor:M_0}, Remark \ref{Rem:M_0}, Lemma \ref{Lem:M_0}, 
and Corollary \ref{pforM}.

The paper is organised as follows. In Section \ref{s:connection} we recall the connection between tuples $(f,g)$ and elliptic surfaces, which will be our main tool.
From a result of Pesenti and Szpiro we derive restrictions on singular 
fibres of such surfaces. These restrictions translate into conditions that 
counterexamples to the Davenport-Stothers inequality must satisfy,
and on the other hand into lower bounds for
$deg(f^3 -g^2)$. Conversely, in Section \ref{s:counterexamples} we
will use base change of elliptic surfaces to construct counterexamples
to $DS(M)$ mod $p$ for many values $M$. (In Section \ref{s:7} we 
completely treat the case $p=7$.) In Section \ref{s:minimising} we 
show that for fixed $p$ and $M$ there is at most one example (up to 
transformation) that reaches the lower bound in 
Theorem \ref{Thm:p} (d). 

In characteristics $2$ and $3$ the connection between $f^3 -g^2$ and
elliptic surfaces does not make sense, but the corresponding question
about fibres of elliptic surfaces does (and is actually more subtle).
We discuss this in Sections \ref{s:3} and \ref{s:2}. Finally,  
Section \ref{s:NS} contains an observation about the N\'eron-Severi group,
which was inspired by the first (conditional) proof of a mod $p$ 
Davenport-Stothers inequality, using a generalised notion of the Artin invariant.

\section{Connection with elliptic surfaces}
\label{s:connection}

The main idea to prove Theorem \ref{Thm:p} is to construct an elliptic
surface $X$ with section over $\PP^1$ which corresponds to the pair of polynomials $(f,g)$. Over a field $k$ of
characteristic $\neq 2,3$, this can be given in Weierstrass form
\begin{eqnarray}\label{eq:Weier}
X = X_{f,g}:\;\; y^2=x^3-3f(t)\,x+2\,g(t)
\end{eqnarray}
with $f, g\in k[t]$. Without loss of generality, we shall assume in the following that $f$ and $g$ have no common zeroes of order at least $4$ resp.~$6$. This implies that the above Weierstrass form (\ref{eq:Weier}) is minimal. Hence the elliptic surface $X$ has singular fibres exactly at the cusps, i.e.~the zeroes of the discriminant
\[
\Delta = 1728 (f^3-g^2)
\]
plus possibly at $\infty$. We employ Kodaira's notation \cite{Kodaira} to 
describe the types of singular fibres. To see how the vanishing orders of $f, g$ and $\Delta$ determine the singular fibre, confer also Tate's exposition in \cite{Tate}. 
If $\deg\ f=2M$ and $\deg\ g=3M$ for some $M\in\N$ and $\deg(f^3 -g^2)<6M$,
then the $j$-invariant of $X$ has a pole of order
$$n=6M-\deg(f^3 -g^2)$$
at $\infty$. More precisely, the singular fibre at $\infty$ has type $I_n$ if $M$ is even
resp.~$I^*_n$ if $M$ is odd. The former is a nodal rational curve for $n=1$ and a  cycle of $n$ lines for $n>1$, corresponding to the extended Dynkin diagram $\tilde A_{n-1}$. The latter corresponds to $\tilde D_{n+4}$.

Conversely, if the fibre at $\infty$ has type $I_n$ or $I_n^*$ with $n>0$, 
then $f$ and $g$ must have degrees $2M$ and $3M$ for some even resp.~odd
$M\in\N$ and $\deg(f^3 -g^2)=6M-n$. Under the condition that the Weierstrass equation is minimal,
the singular fibre at a root of $\Delta$ is additive (i.e.~not of type $I_n$) if and only if this root
is a common zero of $f$ and $g$.

The Euler number of $X$ equals the sum of the Euler numbers of the (singular) fibres. Since these equal the respective vanishing order of $\Delta$, we obtain
$$e(X)=12 \left\lceil\frac M2\right\rceil.$$
There are two further ways to describe the Euler number. First, for an elliptic surface $X$ over $\PP^1$ which is not a product, the first and third Betti numbers vanish: $b_1(X)=b_3(X)=0$. Hence
\begin{eqnarray}\label{eq:e=b_2}
e(X)=b_2(X)+2.
\end{eqnarray}
On the other hand, the canonical divisor $K_X$ is a multiple of the fibre class by \cite[\S 12]{Kodaira}. Hence $K_X^2=0$, and Noether's formula gives for the Euler characteristic $\chi(X)$
\begin{eqnarray}\label{eq:e=chi}
e(X)=12 \,\chi(X).
\end{eqnarray}

If char$(k)=0$, Theorem \ref{Thm:DS} can be proven as follows:
The Shioda-Tate formula \cite[Cor.~5.3]{ShMW} tells us that the N\'eron-Severi group of $X$ is generated by horizontal and vertical divisor classes, i.e.~fibre components and sections. As an $I_n$ fibre has $n$ components, we derive $n\leq \rho(X)-1$. Analogously, it follows for an $I_n^*$ fibre that $n\leq \rho(X)-6$. We shall refer to these estimates as the $\C$-bounds for the singular fibres.

By Lefschetz' theorem, $NS(X)=H^2(X,\Z)\cap H^{1,1}(X)$, so we have for the Picard number
\begin{eqnarray}\label{eq:Lef}
\rho(X)\leq h^{1,1}(X).
\end{eqnarray}
Let $p_g(X)=h^{2,0}(X)$ denote the geometric genus of $X$. Then $\chi(X)=p_g(X)+1$ and $b_2(X)=2\,p_g(X)+h^{1,1}(X)$. From (\ref{eq:e=b_2}) and (\ref{eq:e=chi}), we obtain $h^{1,1}(X)=10 \chi(X)=\frac 56 e(X)$. If the fibre type at $\infty$ is $I_n$ (i.e.~$M$ even), we deduce
\begin{eqnarray}\label{eqn:n}
n\leq \rho(X)-1\leq h^{1,1}(X)-1=\frac 56\, e(X)-1=5\,M-1.
\end{eqnarray}
Hence deg $\Delta=e(X)-n\geq M+1$ proves Theorem \ref{Thm:DS}. The analogous argument applies to type 
$I_n^*$, i.e.~odd $M$.

In positive characteristic, however, we only have Igusa's estimate involving the second Betti number \cite{Igusa}
\begin{eqnarray}\label{eq:b_2}
\rho(X)\leq b_2(X).
\end{eqnarray}
Arguing as in (\ref{eqn:n}) with the weaker bound (\ref{eq:b_2}) and inserting (\ref{eq:e=b_2}) gives the estimate
\begin{eqnarray}\label{eq:weak1}\label{eq:weak}
\text{deg}(f^3-g^2)\geq \begin{cases} 2 & \text{if $M$ is odd,}\\
3 & \text{if $M$ is even,}\end{cases}
\end{eqnarray}
which is independent of $M$. Two of our main aims are to improve this weak estimate (Thm.~\ref{Thm:p} (d)) and to determine when the Davenport-Stothers inequality in its original form holds mod $p$.

The fundamental link between $f^3 -g^2$ and elliptic surfaces is the following

\begin{Lemma}\label{Obs}
Let $f,g$ be as at the beginning of this section. Denote by $X=X_{f,g}/k$ 
the corresponding elliptic surface in Weierstrass form (\ref{eq:Weier}). 
Then the pair $(f,g)$ gives a counterexample to $DS(M)$ mod $p$ if and 
only if the fibre of $X$ at $\infty$ has type at least $I_{5M}^*$, if
$M$ is odd, or $I_{5M}$, if $M$ is even.
\end{Lemma}

\emph{Proof:} Let $I^{(*)}_n$ be the type of the fibre at $\infty$. As discussed earlier, we have $n=6M-\deg(f^3 -g^2)$. Hence the claim of the lemma follows. \qed

Despite its simplicity the following example is, as we shall see later, prototypical for all counterexamples to $DS$ mod $p$.

\begin{Example}[Thm.~{\ref{Thm:p} (a)}]
\label{Ex:deg4}
Let $p>3$ be a prime. We start out with the equation 
$$[t^2 -1]^3 -[t(t^2 -1)]^2 = -(t^2 -1)^2,$$
which of course is not yet a counterexample to $DS$. The corresponding
elliptic surface
$$\widetilde{Y}:\ y^2 = x^3 -3(t^2 -1)x +2t(t^2 -1)$$
has singular fibres of type $I^*_2$ at $\infty$ and $II$ at $\pm 1$. The latter type describes a rational curve with a cusp.
\par
Now let $q$ be any power of $p$ with $q\equiv 1\mod 6$, say $q=6\lambda +1$.
Raising the above equation (over a field of characteristic $p$) to the 
$q$-th power and then dividing by a factor
$(t^2 -1)^{4\lambda}$ resp. $(t^2 -1)^{6\lambda}$, we get
$$f=(t^2 -1)^{2\lambda +1}\ \ and \ \ g=t^q(t^2 -1)\ \ with\ \ 
f^3 -g^2=-(t^2 -1)^2.$$
If $\lambda>1$ this is a counterexample to $DS(M)\ mod\ p$ for 
$M=2\lambda +1=\frac{q+2}{3}$.
In terms of the corresponding elliptic surface we have applied the 
base change $t\mapsto t^q$ to $\widetilde{Y}$ and then minimalised the 
Weierstrass equation of the base-changed surface, which has singular fibres
$I^*_{2q}$, $II$, $II$.
\end{Example}

We note an important property of base change for elliptic surfaces: In the absence of wild ramification (in particular in characteristic $\neq 2,3$), the singular fibres of the base change (after minimalising) can be determined purely in terms of the local ramification data (cf.~\cite{Tate}). For instance, if the local ramification order is denoted by $r\in\N$, then
\begin{eqnarray}\label{eq:I_n}
I_n \to I_{nr} \;\;\;\text{and}\;\;\; I_n^* \to \begin{cases} I_{nr}^*, & \text{if $r$ is odd},\\
I_{nr}, & \text{if $r$ is even}.\end{cases}
\end{eqnarray}
For all other fibre types, the behaviour can be read off from the following diagram:
\begin{eqnarray}\label{eq:add-ram}
II \to \begin{cases} I_0, & \text{if } r\equiv 0 \mod 6,\\
II, & \text{if } r\equiv 1 \mod 6,\\
IV, & \text{if } r\equiv 2 \mod 6,\\
I_0^*, & \text{if } r\equiv 3 \mod 6,\\
IV^*, & \text{if } r\equiv 4 \mod 6,\\
II^*, & \text{if } r\equiv 5 \mod 6,
\end{cases}
\;\;\;\;\;\;
III \to \begin{cases}
I_0, & \text{if } r\equiv 0 \mod 4,\\
III, & \text{if } r\equiv 1 \mod 4,\\
I_0^*, & \text{if } r\equiv 2 \mod 4,\\
III^*, & \text{if } r\equiv 3 \mod 4.
\end{cases}
\end{eqnarray}
Here we employ the convention that $I_0$ describes a smooth curve of genus one.

For the next examples, and actually 
throughout the paper, we fix another rational elliptic surface which we denote by $Y$. It is given in Weierstrass form
\begin{eqnarray}\label{eq:Y}
Y:\;\;y^2 = x^3 - 3\, t^3(t-1) x + 2\, t^5(t-1).
\end{eqnarray}
This has discriminant $\Delta=-1728\, t^9 (t-1)^2$ and 
singular fibres of types $III^*$ at $0$, $II$ at $1$ and $I_1$ at 
$\infty$. For any $p>3$,
$Y$ has good reduction mod $p$. In the following, we always start
with $Y$ over some fixed field $\F_p$.
\par
Note that the surface $\widetilde{Y}$ in the previous example can be
obtained from $Y$ as follows: First apply the degree $2$ base change 
$t\mapsto t^2$. The resulting surface has the following bad fibres:
$I_2$ at $\infty$, $I^*_0$ at $0$, $II$ at $1$ and $-1$. Now apply
the quadratic twist that ramifies exactly at $\infty$ and $0$ (i.e.~remove the common factors $t^2, t^3$ of $f, g$).

\begin{Example}\label{Ex:25}
Let $p=5$. Let $X$ be the pull-back of $Y/\F_p$ by the purely
inseparable base change $t\mapsto t^{25}$. Then $X$ has singular
fibres $I_{25}, III^*, II$. Thus it produces a counterexample to
$DS(5)$ mod $p$: Eliminating the common factor $t$ (to the order
$2$ resp.~$3$), we obtain
\[
[t^{14}(t-1)]^2-[t(t-1)^9]^3\equiv t^3(t-1)^2 \mod 5.
\]
In terms of elliptic surfaces, the elimination of $t^2$ resp. $t^3$ again
corresponds to a quadratic twist that ramifies exactly at 
$0$ and $\infty$; so it changes the configuration to $I^*_{25}$, $III$, $II$.
\end{Example}

\begin{Example}\label{Ex:31}
Let $p=31$. Let $X$ be the pull-back of $Y/\F_p$ by the purely
inseparable base change $t\mapsto t^{31}$. Then $X$ has singular
fibres $I_{31}, III, II$. The counterexample to $DS(6)$ mod $p$
is:
\[
[t^{17}(t-1)]^2-[t(t-1)^{11}]^3\equiv t^3(t-1)^2 \mod 31.
\]
\end{Example}

Note that we can apply the above constructions to any prime-power
$q=p^r\equiv 1\mod 6$ (and with other singular fibres resulting also to $q=p^r\equiv 5\mod 6$). For fixed $p$, this gives infinitely many
counterexamples to $DS(M)$ mod $p$ for some wide-spread $M$. To prove 
Thm.~\ref{Thm:p} (f), it thus remains to fill in the gaps in between 
these $M$; this will be done in Section \ref{s:counterexamples}.

The above examples are all unirational, as C.~Schoen pointed out. We will later see that any counterexample with deg$(f^3 -g^2)\le 6$ has to be unirational (Prop.~\ref{Prop:unique}). What makes the above constructions work are two properties: On the one hand, the number of cusps stays unchanged under the purely inseparable base change. On the other hand, additive fibres of type $\neq I_n^* (n>1)$ are interchanged according to (\ref{eq:add-ram}). Hence their contribution to the Euler number is always small (since $\Delta$ vanishes at most to order $10$ at these fibres). Hence, the $I_n^{(*)}$ fibre at $\infty$ becomes more and more dominant under (\ref{eq:I_n}). This dominance is unbroken if we compose with suitable separable base changes. This will be exhibited in Sec.~\ref{s:counterexamples}.

\section{Conditions on counterexamples}
\label{s:conditions}

As we have seen in the previous section, giving lower bounds for 
$deg(f^3 -g^2)$ corresponds in terms of the elliptic surface $X_{f,g}$ 
and the fibre $I^{(*)}_n$ at $\infty$ to giving upper bounds for $n$ by the Euler number $e(X_{f,g})$. In order to do so, we use a result of 
Pesenti-Szpiro \cite{P-S}, which we recall in the following in slightly 
simplified setting.

Let $X$ be a non-isotrivial elliptic surface over $\PP^1$. Denote by $p^d$ 
the degree of inseparability of the map $X\to\PP^1$. That is, we have a 
commutative diagram 

$$
\begin{array}{ccc}
X & \to & S(X)\\
&&\\
\downarrow && \downarrow\\
&&\\
\PP^1 & \to & \mathbb{P}^1
\end{array}
$$
where the elliptic surface $S(X)\to\PP^1$ is separable and the map $\PP^1\to\PP^1$ is purely inseparable of degree $p^d$. Let $\NN=\NN_X$ be the conductor of $X$, so that at $v\in\PP^1$

\[
\text{ord}_v(\NN) = \begin{cases} \;\;\;0,\\ \;\;\;1,\\ 2+\delta_v,\end{cases} \text{if the fibre $X_v$ is}\; 
\begin{cases} 
\text{smooth} \;\,(I_0),\\ 
\text{multiplicative}\;\, (I_n, n>0),\\ 
\text{additive} \;\,(\neq I_n).
\end{cases}
\]
Here $\delta_v$ measures the wild ramification at $v$. In particular, $\delta_v=0$ if char$(k)\neq 2,3$. Note also that $\NN | \Delta$ and that $\NN_X=\NN_{S(X)}$ by inseparability.

\begin{Theorem}[Pesenti-Szpiro {\cite[Thm.~0.1]{P-S}}]
\label{Thm:P-S}
With the above notation, we have
\[
e(X)\leq 6 \,p^d (\text{deg}\, \NN -2).
\]
\end{Theorem}

\begin{Corollary}\label{Cor:P-S}
Let $X$ be a separable elliptic surface over $\PP^1$. Then its fibres obey the $\C$-bounds.
\end{Corollary}

\emph{Proof:} Let $X$ have a fibre of type $I_n$. Then 
\[
\text{deg}\, \NN \leq 1 + (e(X) - n)
\]
since $\NN | \Delta$. By Thm.~\ref{Thm:P-S},
\[
e(X) \leq 6 (e(X) - n - 1).
\]

This simplifies to $n\leq	\frac 56 e(X) - \frac 16 < \frac 56 e(X)$ as claimed. In the case of a fibre of type $I_m^*$, deg $\NN\leq 2 + (e(X) - m - 6)$, and the same argumentation applies. \qed

\begin{Corollary}\label{Cor:sep}
Any counterexample to $DS$ occurs via purely inseparable base change.
\end{Corollary}

\begin{Remark}[Thm.~{\ref{Thm:p} (b)}]
\label{Rem:multiplicity}
Any factor of $f^3 -g^2$ that does not divide $f$ and $g$ corresponds to 
a fibre $I_n$ of $X_{f,g}$. If $(f,g)$ is a counterexample to $DS$, then
the corollary immediately implies that $n$, which is also the multiplicity 
of this factor in $f^3 -g^2$, is divisible by $p$. On the other hand, every 
common factor of $f$ and $g$ has of course at least multiplicity $2$ in
$f^3 -g^2$. Thus, in any counterexample, $f^3 -g^2$ has no factors with 
multiplicity $1$ at all.
\end{Remark}

As a further result we note that the Davenport-Stothers inequality stays 
valid  in characteristic $p$ for polynomials which are coprime or have only one common zero:

\begin{Proposition}[Thm.~{\ref{Thm:p} (c)}]
\label{Lem:coprime}\label{Prop:s-s}
Let $p>3$. Let $f,g \in k[t]$ with deg $f=2M$, deg $g=3M$ and $f^3\neq g^2$.
Then
\[ 
\# \{s\in k; f(s)=g(s)=0\}\leq 1  \Rightarrow \text{deg}(f^3-g^2)\geq M+1.
\]
\end{Proposition}

\emph{Proof:} By assumption, the corresponding elliptic surface $X_{f,g}$ has only multiplicative fibres away from $\infty$, except for possibly one additive fibre. Hence the proposition follows from the following general fact which can already be derived from a theorem of Szpiro \cite[Thm.~1]{Szpiro} preceeding his later work with Pesenti \cite{P-S}. \qed

\begin{Lemma}
\label{Lem:s-s}
Let char$(k)\neq 2, 3$. Then any surface over $\PP^1$ with only $I_n^{(*)}$ fibres plus possibly one other singular fibre obeys the $\C$-bounds.
\end{Lemma}

\emph{Proof:}  We prove the corollary in three steps, assuming that there is an elliptic surface $X$ as above which disobeys the $\C$-bounds.

\emph{1. reduction step: elimination of  $I_n^*$ fibres}

Let $X'$ arise from a quadratic twist which only ramifies over $I_n^*$ fibres. The twist decreases the Euler number, so $X'$ also disobeys the $\C$-bounds. Hence we can assume without loss of generality that $X$ has at most one additive fibre other than the special fibre at $\infty$.

\emph{2. reduction step: elimination of additive fibres}

Let $X'$ be a base change of degree $12$ of $X$ which only ramifies above $\infty$ and above the single additive fibre (or above an arbitrary point, if there is no additive fibre away from $\infty$). By construction, $X'$ is semi-stable, i.e.~there are no additive fibres. Since $e(X')\leq 12\, e(X)$, it is easily checked that $X'$ disobeys the $\C$-bounds. Hence we can assume that $X$ itself is semi-stable.

\emph{3. proof for semi-stable fibrations}

If $X$ is separable, we are done. Otherwise consider the separable elliptic surface $S(X)$ as before. Then $e(X)=p^d\, e(S(X))$. By Cor.~\ref{Cor:P-S}, a fibre of type $I_n$ on $S(X)$ satisfies $n<\frac 56 e(S(X))$. Hence, the resulting fibre of type $I_{p^dn}$ on $X$ satisfies $p^d\, n<\frac 56\, p^d\, e(S(X)) = \frac 56\, e(X)$. \qed

\begin{Remark}\label{Rem:s-s-2,3}
In characteristic $2$ and $3$, Lem.~\ref{Lem:s-s} holds in the absence of 
wild ramification. In particular, this applies to semi-stable fibrations.
\end{Remark}

\begin{Remark}
Using the techniques of the proof of Lemma \ref{Lem:s-s}, one can strengthen Prop.~\ref{Lem:coprime} as follows: Let $p>3$. Let $f,g \in k[t]$ with deg $f=2M$, deg $g=3M$ and $f^3\neq g^2$.
Then
\[
\{s\in k; f(s)=g(s)=0\}=\{s_0\} \Rightarrow \text{deg}(f^3-g^2)\geq M+\text{ord}_{s_0}(f^3, g^2).
\]
\end{Remark}

\section{Mod $p$ Davenport-Stothers inequalities}
\label{s:inequality}

In this section we discuss which bounds for deg$(f^3-g^2)$ do exist
in characteristic $p>3$ (instead of the Davenport-Stothers inequality from Thm.~\ref{Thm:DS}). In terms of the corresponding elliptic
surface $X$, we a priori only have Igusa's estimate (\ref{eq:b_2}). In
case $f^3\neq g^2$, this estimate directly implies the bound (\ref{eq:weak})
\begin{eqnarray*}
\text{deg}(f^3-g^2)\geq \begin{cases} 2 & \text{if $M$ is odd,}\\
3 & \text{if $M$ is even.}\end{cases}
\end{eqnarray*}
For $X$ rational, i.e.~$M=1,2$, this bound is sharp. However, for
$M>2$, an improvement is given by the following proposition. 

\begin{Proposition}[Thm.~{\ref{Thm:p} (d)}]
\label{Prop:4,5}
Let $p>3$, $M>2$ and $f,g \in k[t]$ with deg $f=2M$, deg
$g=3M$. Then
\[ 
f^3\neq g^2 \Leftrightarrow \text{deg}(f^3-g^2)\geq \begin{cases} 4 & \text{if $M$ is odd,}\\
5 & \text{if $M$ is even, } p\equiv 7\text{ mod } 12\;\; \text{ or } M=4,\\
6 & \text{if $M$ is even, } M>4,\ p\not\equiv 7 \text{ mod } 12.
\end{cases}
\]
\end{Proposition}

\emph{Proof:} Consider the associated elliptic surface $X=X_{f,g}$, so $e(X)\leq 12\lceil\frac M2\rceil$. The argumentation of Lem.~\ref{Obs} shows that Prop.~\ref{Prop:4,5} is equivalent to the following lemma applied to $X$. \qed

\begin{Lemma}\label{Lem:4,5}
Let $p>3$ and $X$ an elliptic surface over $\PP^1$ in characteristic $p$ with Euler number $e$. If $X$ has a fibre of type $I_{e-3}, I_{e-4}, I_{e-8}^*$ or $I_{e-9}^*$, then $X$ is rational (i.e.~$e=12$). If $X$ has a fibre of type $I_{e-5}$, then $X$ is rational or K3 or $p\equiv 7\mod 12$.
\end{Lemma}

The case of $e=24$ (thus $M=3,4$) has been proven in \cite{S-max} for all $p$. The proof therein used Artin's theory of supersingular K3
surfaces \cite{A} and in particular the \emph{Artin invariant}.
Here, we could combine the techniques from \cite{S-max} with Thm.~\ref{Thm:NS} to prove Lem.~\ref{Lem:4,5} for types $I_{e-3}$ and $I_{e-8}^*$. However,  Thm.~\ref{Thm:P-S} gives us complete control of the situation.

\emph{Proof of Lemma~\ref{Lem:4,5}:} 
Let $X$ be as in the lemma and exclude fibre type $I_{e-5}$ for the time being. Then deg $\NN_X\leq 5$. By Thm.~\ref{Thm:P-S}, the separable elliptic surface $S(X)$ is rational. So if $X$ is separable, we are done. 

Assume $X\neq S(X)$. Then $X$ gives a counterexample to $DS$ mod $p$. Consider the singular fibres other than the special one from the lemma. Then $X$ cannot have another fibre of type $I_n^{(*)}$, since $n\leq 4<p$. Hence, all other singular fibres are additive. By Prop.~\ref{Lem:coprime}, $X$ has at least three singular fibres. Hence their configuration is $I_{e-8}, II, II$. Thus $S(X)$ can only have singular fibres $I_8, II, II$ or $I_4, IV, IV$. Both configurations are impossible. To see this, apply a quadratic twist which ramifies at the two additive fibres. The resulting configuration $I_8, IV^*, IV^*$
resp. $I_4, II^*, II^*$ contradicts Theorem \ref{Thm:P-S}.

Finally, assume that $X$ has a fibre of type $I_{e-5}$, but is not rational or K3. In particular, $X$ is not separable by Cor.~\ref{Cor:sep}. Thus $X$ has no other multiplicative singular fibre. By Prop.~\ref{Lem:coprime}, there are two further cusps, so the configuration is $I_{e-5}, III, II$. This implies that deg $\NN=5$, so $S(X)$ is rational. It is immediate that $S(X)$ can only have configuration $I_7, III, II$ or $I_1, III^*, II$. However, the former configuration is ruled out for a separable elliptic surface by quadratic twisting as above. We obtain for the degree of inseparability
\[
p^d = e-5 \equiv 7\mod 12.
\]
This proves the claim. \qed

\begin{Remark}\label{Rem:sharp}
In general, the bounds of Prop.~\ref{Prop:4,5} are sharp. To see this for the second, apply any purely inseparable base change to $Y$ as in Ex.~\ref{Ex:31}. For the first and third, consider any prime power $q=p^r\equiv 1\mod 12$. Then pull-back from $\widetilde Y$ via $t\mapsto t^q$ as in Ex.~\ref{Ex:deg4} or compose with the base change $\pi_3: t\mapsto 2\,t^3-1$.
\end{Remark}

Combining 
Lem.~\ref{Obs}, Cor.~\ref{Cor:sep} and Prop.~\ref{Prop:4,5}
we obtain the following sufficient, but in general not necessary
criterion.

\begin{Criterion}[Thm.~{\ref{Thm:p} (e)}]
\label{Crit}
$DS(M)$ holds mod $p$ if there is no $n\in\N$ such that
\[
5M\leq np\leq 6M-\begin{cases} 4, & \text{if $M$ is odd,}\\
5, & \text{if $M$ is even},\, p\equiv 7\text{ mod } 12\\
6, & \text{if $M$ is even},\, p\not\equiv 7 \text{ mod } 12. \end{cases}
\]
\end{Criterion}

\section{Counterexamples for almost all $M$}
\label{s:counterexamples}

The main idea of this section is the
following observation: Let $X$ be a counterexample for $DS(M)$
mod $p$. Often $X$ will also give rise to a counterexample for
$DS(M+1)$ mod $p$. This can be achieved by quadratic twisting, i.e.~adding a fibre of type
$I_0^*$ or transferring the $*$. In terms of $f$ and $g$, this
corresponds to adding a common (linear) factor:
\[
f \mapsto \alpha^2 f,\;\;\; g\mapsto \alpha^3 g.
\]

We shall now make this explicit: Fix the characteristic $p>3$ and
let $q=p$ or $p^2$ such that $q\equiv 1 \mod 12$ and $q>13$. Let
$X$ denote the pull-back of $Y$ via $t\mapsto t^q$. This has
singular fibres $I_q, III^*$ and $II$ (as in Ex.~\ref{Ex:25}). It
is easy to see that we obtain a counterexample to $DS(M)$ mod $p$
with odd $M=\frac{q+5}6$. This suffices for $DS(M+1)$ mod $p$ as
explained above if $q>49$. Then we shall ask whether it also suffices for $DS(M+2)$ mod $p$, or search for other adequate surfaces.

In the following, we will apply the cyclic base changes
\[
\phi_d: t \mapsto t^d
\]
to our fixed surface $X$. Denote the resulting elliptic surface by $X_d$.  By (\ref{eq:I_n}), (\ref{eq:add-ram}), the singular fibres of $X_d$ are as follows:
\begin{eqnarray}\label{eq:X_d}
\begin{array}{ll}
d\equiv 0\mod 4: & I_{qd},  d ~II,\\
d\equiv 1\mod 4: & I_{qd}, III^*, d~ II,\\
d\equiv 2\mod 4: & I_{qd}, I_0^*, d ~II,\\
d\equiv 3\mod 4: & I_{qd}, III, d~ II.\\
\end{array}
\end{eqnarray}

We can read off from the
singular fibres, that $X_d$
contradicts $DS(M_d)$ mod $p$ for a certain $M_d$ which
asymptotically equals $\frac{e(X_d)}{6} \sim \frac{d(q+2)}6$. To prove 
the first statement of Theorem \ref{Thm:p} (f), we only need the following

\begin{Lemma}\label{Lem:d_0}
There is some $d_0\in\N$, such that for all $d\geq d_0$, $X_d$
gives counterexamples to $DS(M)$ mod $p$ for all
$M=M_d,\hdots,M_{d+1}-1$
\end{Lemma}

In other words, for $d\geq d_0$, $X_d$ fills all the gaps between
$M_d$ and $M_{d+1}$ with counterexamples. Hence we can choose
$M_0=M_{d_0}$ to deduce Theorem \ref{Thm:p} from the above lemma.

\emph{Proof of Lemma \ref{Lem:d_0}}: From Lem.~\ref{Obs}, it is easy to deduce that
the lemma requires
\[
qd>5(M_{d+1}-1)-1.
\]
Reading off the Euler number of $X_d$ from the singular fibres in each case of (\ref{eq:X_d}), we immediately obtain the bound
\[
M_d\leq \frac 16 (d(q+2)+3).
\]
Hence it suffices to check that
\[
dq>5\left(\frac 16 ((d+1)(q+2)+3)-1\right)-1.
\]
Equivalently
\[
d(q-10)>5q-11,
\]
so we need
\begin{eqnarray}\label{eq:d=6}
d>5+\frac{39}{q-10}.
\end{eqnarray}
Lemma \ref{Lem:d_0} follows. \qed

\begin{Corollary}[Thm.~{\ref{Thm:p}} (f)]\label{Cor:M_0}
Let $q=p^r>29$ and $M_0=q+2$. Then $DS(M)$ mod $p$ does not hold if $M\geq M_0$.
\end{Corollary}

\emph{Proof:} By inequality (\ref{eq:d=6}) in the proof of Lemma \ref{Lem:d_0}, we can always choose $d_0=6$ if $q>49$ and $q\equiv 1\mod 12$. The same arguments with fibre types $III$ and $III^*$ interchanged apply to $q\equiv 7\mod 12$. The corresponding elliptic surface $X_6$ has the following configuration and Euler number:
\begin{eqnarray*}
I_{6q}^*, 6 \,II,\; & e(X_6)=6 (q+3), & M_6=q+2,\;\;\;\text{if } q\equiv 1\mod 6,\\
I_{6q}^*, 6 \,IV, & e(X_6)=6 (q+5), & M_6=q+4,\;\;\;\text{if } q\equiv 5\mod 6.
\end{eqnarray*}

If $q\equiv 1\mod 6$ and $q>49$, Corollary \ref{Cor:M_0} thus follows. For $q=31, 37, 43, 49$, however, we only have $d_0=7$. It remains to investigate the gap between $X_6$ and $X_7$. 

Note that for $q=49$, there is in fact no gap, since $\phi_7$ becomes purely inseparable mod $7$: This implies that $X_7$ has configuration $I_{343}, III, II$, so $M_7=\frac{348}6=58$. On the other hand, the $I_{6q}^*$ fibre of $X_6$ suffices up to $M=\lfloor\frac{6q}{5}\rfloor=58$ by Lem.~\ref{Obs}. Hence there is no gap.

We shall treat one further case in detail.
Let $q=31$. Then $X_6$ with fibre $I_{186}^*$ provides counterexamples for $M=33,\hdots,37$.  Meanwhile, $X_7$ has configuration $I_{217}^*, III, 7\,II$, so $M_7=39$. To obtain a counterexample for $M=38$, we use a base change which unlike $\phi_7$, is not cyclic. Instead, let $\tilde\pi_7$ have ramification indices $7$ at $\infty$, $(2,2,2,1)$ at
$0$ and $(3,2,1,1)$ at $1$. The resulting elliptic surface has configuration $I_{7q}, IV, III, 2\,II$ and Euler number $228$, so it in fact gives a counterexample for $M=\frac{228}6=38$.

Over $\overline{\Q}$, there are three such base changes. These are conjugate 
in the extension
$\Q(x^3-x^2+5x+1)=\Q(\alpha^3+7\alpha^2+91\alpha+315)$ and can be defined 
as follows:
\[
\tilde\pi_7: t\mapsto \left(1+\frac 16\alpha\, t+\frac1{24}\alpha^2\,t^2-\frac 7{120}(\alpha^2+6\alpha+60)\,t^3\right)^2\;\left(1-\frac 13\alpha\, t\right)
\]
All these base changes have good reduction outside the primes above $\{2,3,7\}$.
Applying $\tilde\pi_7$ to $X/\F_{31}$ gives the required counterexample for $M=38$. 

It is easily checked that $\tilde\pi_7$ also suffices to fill the gap between $X_6$ and $X_7$ for $q=37, 43$. This concludes the proof of Corollary \ref{Cor:M_0} for $q\equiv 1\mod 6$.

For $q\equiv 5\mod 6$, we furthermore need the following base changes:
\begin{itemize}
\item $\pi_8$ with ramification indices $8, (2,2,2,2),
(3,2,1,1,1)$. This can be obtained from $\phi_2$ as above and $\pi_4$ with ramification indices $4, (3,1), (2,1,1)$ as in
\cite[\S 3]{S-RM} by way of reduction.  
\item $\pi_6: t\mapsto (2t^3-1)^2$ with
ramification indices $6, (2,2,2), (3,1,1,1)$ (which factors through $\pi_3$ from Rem.~\ref{Rem:sharp}).\qed
\end{itemize}

\begin{Remark}[Thm.~{\ref{Thm:p}} (f)]\label{Rem:M_0}
For $q\equiv 1\mod 6$, the above bound can be improved to $M_0=\frac{5q+7}{6}$ by considering the pull-back of $X$ via $\pi_6$ as above and via
$\pi_H$ with ramification indices $5, (3,1,1), (2,2,1)$ from \cite[\S 7]{S-RM}. Moreover, this bound also holds for $q=
25$. We shall see in Lemma \ref{Lem:M_0} that these values are sharp.
\end{Remark}

\begin{Remark}
We can always find counterexamples for $DS$ mod $p$ with
coefficients in $\F_p$. 
\end{Remark}

\emph{Proof:} The only base change which is not a priori defined over $\F_p$, is $\tilde\pi_7$, since it does not arise from $\Q$ by reduction. Note that the Galois extension $\Q(x^3-x^2+5x+1)/\Q$ has $\Q(\sqrt{-3})$ as
intermediate field. Hence, the polynomial $x^3-x^2+5x+1$ factors into a linear and a quadratic factor mod $p$, if $p$ is inert in $\Q(\sqrt{-3})$. In other words, we find $\tilde\pi_7$ over $\F_p$ if $p\equiv 5\mod 6$.

For $p\equiv 1\mod 6$, it is useful to work with a base change $\pi_7$ with ramification indices $7$ at $\infty$, $(2,2,2,1)$ at
$0$ and $(3,1,1,1,1)$ at $1$ instead. The resulting elliptic surface has different configuration $I_{7q}, III, 4\,II$, but the Euler number stays the same, so it also gives a counterexample for $M=\frac{228}6=38$.

Over $\C$, such base changes come in a one parameter-family which specialises to $\tilde\pi_7$. We computed one parametrisation of the family explicitly as
\[
\pi_\lambda: t\mapsto \left(1+\frac 12\,t+\frac 38\,t^2-\lambda\,t^3\right)^2\;(1-t)
\]
with $\lambda\neq 0, \frac{15}8$. Specialising, for instance, at $\lambda=-\frac3{16}$, we obtain a base change $\pi_7$ which lives over $\Q$ with a further rational cusp at $t=-2$  and good reduction outside $\{2,3\}$. \qed

Criterion \ref{Crit} helps us establish the optimality 
of $M_0$.

\begin{Lemma}[Thm.~{\ref{Thm:p} (f)}]
\label{Lem:M_0}
If $p\equiv 1\mod 6$, then $DS(\frac{5p+1}{6})$ holds mod $p$. 
If $p\equiv 5\mod 6$, then $DS(p+1)$ holds mod $p$. 
In particular, for $p>29$ the value $M_0$ defined in Cor.~\ref{Cor:M_0} 
resp.~Rem.~\ref{Rem:M_0} is optimal.
\end{Lemma}

\emph{Proof:} For $p\equiv 1\mod 6$ this follows immediately from
Criterion \ref{Crit}. 

If $p\equiv 5\mod 6$, again by Criterion \ref{Crit} the corresponding 
elliptic surface $X$ must have a fibre of type $I_{6p}$. Furthermore, the configurations of $X$ and $S(X)$ have to be among the following:
$$
\begin{array}{ccccccc}
X:&& I_{6p}\;\;\; \text{ and}&
II, II, II&
\;\;IV, II&
\;\;III, III&
\;\;I_0^*,\\
S(X):&& I_{6}\;\;\; \text{ and}&
II^*, II^*, II^* &
\;\;II^*, IV^*&
\;\;III^{(*)}, III^{(*)}&
\;\;I_0^*.
\end{array}
$$
All but the reduced case of the third configuration of $S(X)$ give direct contradictions by Thm.~\ref{Thm:P-S}. The same argument applies to the quadratic twist of the reduced case of the third configuration which ramifies above the additive fibres. \qed

\section{$DS$ mod $7$}
\label{s:7}

This section gives a complete study of $DS$ mod $7$. By
Remark \ref{Rem:M_0}, applied to $q=49$, $DS(M)$ does not hold mod
$7$ for $M\geq 42$. The following table collects all further
counterexamples.

The first column lists the degree $M$. It is followed by the base
changes used to construct a counterexample from the rational
elliptic surface $Y/\F_7$. To compose those base changes which are
separated by a comma, we need an automorphism of $\PP^1$
translating some cusps. We omit the details, since these can be
extracted from the resulting fibre configuration. This is
displayed in the third column.

\begin{table}[ht!]
\begin{tabular}{cc}

{\begin{tabular}{|c||c|c|} \hline
$M$ & base changes & configuration\\
\hline\hline $\vdots$ &&\\
42 & $\phi_{49}\circ\pi_5$ & [245, III, 2 II] \\
\hline
41 &  &  \\
40 & $\phi_7\circ\pi_5, \phi_6$ & [210*, 12 II] \\ \hline
39 &  &  \\
$\vdots$ &  &  \\
34 & $\phi_{49}\circ\phi_4$ & [196, 4 II] \\ \hline
33 &  &  \\
32 &  &  \\
31 & $\phi_7\circ\pi_6, \phi_4$ & [168*, 9 II] \\ \hline 30 &
$\phi_7\circ\phi_2, \phi_{11}$ & [154, IV, 11 II] \\ \hline
29 &  &  \\
$\vdots$ &  &  \\
26 & $\phi_{49}\circ\phi_3$ & [147, III, 3 II] \\ \hline
25 &  &  \\
24 &  &  \\
23 & $\phi_7\circ\phi_2, \phi_3, \phi_3$ & [126*, 6 II] \\ \hline
\end{tabular}}

&

{\begin{tabular}{|c||c|c|} \hline
$M$ & base changes & configuration\\
\hline\hline
22 &  &  \\
21 & $\phi_7\circ\phi_2, \pi_4, \phi_2$ & [112*, 7 II] \\ \hline
20 & $\phi_7\circ\pi_5, \phi_3$ & [105, III, 6 II] \\ \hline
19 &  &  \\
18 &  &  \\
17 & $\phi_{49}\circ\phi_2$ & [98*, 2 II] \\ \hline 16 &
$\phi_7\circ\phi_2, \phi_6$ & [84, 6 II] \\ \hline 
15 & $\phi_7\circ\pi_{11}$ & [77*, III, 5 II] \\
\hline
14 &  &  \\
13 & $\phi_7\circ\phi_2, \tilde\pi_5$ & [70*, 4II] \\ \hline 
12 & $\phi_7\circ\pi_9$ & [63, III, 3 II] \\ \hline 11 & $\phi_7\circ\phi_2, \phi_4$ & [56*, 5 II] \\
\hline 10 & - &  \\ \hline 9 & $\phi_{49}$ & [49*, III, II] \\
\hline 8 & $\phi_7\circ\pi_6$ & [42, 3 II] \\ \hline 7 &
$\phi_7\circ\pi_5$ & [35*, III, 2 II] \\ \hline 6  & -  &  \\
\hline
5  & - &  \\
\hline
\end{tabular}}
\end{tabular}
\caption{Counterexamples to $DS(M)$ mod $7$}
\end{table}

In the table, there are three base changes which have not been introduced yet. One of them,
\[
\tilde\pi_5: t \mapsto  t^3\,(6t^2-15t+10),
\]
does already have local ramification indices $5, (3,1,1), (3,1,1)$ when regarded in zero characteristic. The other two base changes were only computed modulo $7$. We list them with the ramification indices:
$$
\begin{array}{rcllll}
\pi_9: t & \mapsto & (t^3+5 t^2+6 t+4)^2 (t+4)^2 t & \;\;9, (2,2,2,2,1), (3,3,1,1,1)\\

\pi_{11}: t & \mapsto & -t(t^4+6t^3+5t^2+4)^2(t+6)^2 & \;\;11, (2,2,2,2,2,1), (3,3,1,1,1,1,1)
\end{array}
$$

As a consequence, all base changes involved can be defined over $\F_7$. 
The same can be checked for the compositions. The remaining gaps in the table are clarified by the following

\begin{Lemma}\label{Lem:7}
$DS(M)$ holds mod $7$ if and only if $M\in\{1,2,3,4,5,6,10\}$.
\end{Lemma}

\emph{Proof:} The table shows that $DS$ mod $7$ does not hold for any 
other $M$. The other implication of Lemma \ref{Lem:7} can easily be verified 
using Criterion~\ref{Crit}.
\qed

\begin{Remark}
With similar methods one could for any given $p$ find all values
$M$ such that $DS(M)$ holds mod $p$. For example, $DS(M)$ holds mod $31$
if and only if 
\[
M\in\{1,2,3,4,5,7,8,9,10,13,14,15,16,19,20,21,25,26\}.
\]
\end{Remark}

Finally, we consider the converse problem, starting with a fixed $M$.
Then it is clear from Shioda's original observation (cf.~Criterion \ref{Crit}) that $DS(M)$ holds mod $p$ 
for almost all $p$. For example $DS(5)$ mod $p$ holds for all $p>5, p\neq 13$.
Nevertheless:

\begin{Corollary}[Thm.~{\ref{Thm:p} (g)}]
\label{pforM}
Let $M>4$. Then there is a $p>3$ such that $DS(M)$ does not hold
mod $p$. Moreover, a counterexample can be given over $\F_p$.
\end{Corollary}

\emph{Proof:} For $M\neq 5,6,10$, this follows from the table for $p=7$. 
For the remaining values of $M$, consider Examples \ref{Ex:deg4} and \ref{Ex:31} as well as the 
base change of $Y/\F_{53}$ via $\phi_{53}$. \qed

\section{Multiplicative fibres in counterexamples}
\label{s:multi}

All counterexamples which we have constructed so far, have exactly one fibre of type $I_n^{(*)}$ with $n>0$. The only reason for this restriction is that our techniques are most effective in this setting. In general, the following result is easily derived from Lemma \ref{Obs}.

\begin{Lemma}
\label{Lem:multfibres}
Suppose that $X$ has a fibre of type $I^{(*)}_n$ at $\infty$.
A necessary and sufficient condition that $X$ gives rise to a counterexample
to $DS$ mod $p$ after a purely inseparable base change of sufficiently large degree is
$$n>5\sum_{\nu\in T} m_{\nu}$$
where $T$ is the set of fibres $I^{(*)}_{m_{\nu}}$ outside $\infty$.
\end{Lemma}

\begin{Example}
Consider the rational elliptic surfaces with configuration [7, 1, 2 II] and [6, 1, III, II] from \cite{Herf}. Both have good reduction at (the primes above) $p\neq 2,3,7$. Then apply a purely inseparable base change.
\end{Example}

One can easily determine which degrees suffice for the above examples. In general, we combine Remark \ref{Rem:multiplicity} and Proposition 
\ref{Lem:coprime} to deduce that a counterexample to $DS(M)$ mod $p$ with $p> M-4$ cannot
have a multiplicative fibre outside $\infty$.

Now we show that (at least) at the price of having many additive
fibres, we can construct counterexamples to $DS$ mod $p$ with as
many multiplicative fibres as we want, and we can even choose where
to locate these multiplicative fibres.

\begin{Example}
Choose a finite set $V\subseteq k^*$ of places $\nu$
together with positive integers $m_{\nu}$ and $n$ so that
$n>5\sum_{\nu\in V}m_\nu$. Let $d=n+\sum_{\nu\in V}m_\nu$. Define the base change $\pi$ of degree $d$ by
\[
\pi: t\mapsto \dfrac{t^d}{\prod_{\nu\in V} (t-\nu)^{m_\nu}}.
\]
Apply $\pi$ to the elliptic surface $Y$ from (\ref{eq:Y}). This base change results in a fibre $I_n$ at
$\infty$ and fibres $I_{m_{\nu}}$ at $\nu$ plus generally $d$ fibres of type $II$ and the special fibre at $0$ whose type depends on the residue class of $d$ mod $4$. By construction, the pull-back satisfies the condition of Lem.~\ref{Lem:multfibres}.
\end{Example}

In deriving an actual counterexample from such a base change, the large number of additive fibres causes a large degree of inseparability. Hence it would be interesting to bound the number or types of additive fibres. Recall that there are at least two additive fibres (outside $\infty$) by Prop.~\ref{Prop:s-s}. Here we shall only treat the first few cases. For this purpose, we introduce the following notation for an elliptic surface $X$ over $\PP^1$:

\begin{center}
\begin{tabular}{ll}
$F_v$ & fibre above $v\in\PP^1$\\
$m_v$ & number of components of $F_v$\\
$U$ & set of multiplicative cusps other than $\infty$\\
$S$ & set of additive cusps other than $\infty$
\end{tabular}
\end{center}

\begin{Lemma}
If $X$ is separable and every fibre outside $\infty$ is reduced,
then the condition in Lem.~\ref{Lem:multfibres} implies
\[
6+6\sum_{\nu\in U} (m_\nu-1)<\sum_{\nu\in S} (5-m_\nu).
\] 
\end{Lemma}

\emph{Proof}:  We shall only treat one case in detail, say $F_\infty$ 
reduced, $\# S$ even. Consider the quadratic twist $X^*$ of $X$ above $S$. 
Then Theorem \ref{Thm:P-S} gives
\[
e(X^*)=m_\infty + \sum_{\nu\in U} m_\nu +\sum_{\nu\in S} (m_\nu+7) \leq 6\, (\# U+2\,\# S-1).
\]
Inserting the estimate for $m_\infty$ from Lem.~\ref{Lem:multfibres}, we obtain
\[
6 \sum_{\nu\in U} m_\nu +\sum_{\nu\in S} (m_\nu+7) < 6 \,(\# U+2\,\# S-1).
\]
The claim follows immediately. The other cases can be proved analogously, 
twisting over $S\cup\{\infty\}$ if $\# S$ is odd.
\qed

\begin{Corollary}
Under the above hypotheses, $X$ cannot have
the following configurations of additive fibres (other than at $\infty$):
\[
2 IV\;\;-\;\; IV, III\;\;-\;\;IV, II\;\;-\;\;2 III\;\;-\;\; 3 IV.
\]
\end{Corollary}

\begin{Corollary}
Let $X$ be a counterexample to $DS$ mod $p$. Then the additive configurations 
$[II, IV]$ and $[2 III]$ and any of their quadratic twists are impossible. The additive configurations $[3 IV], [2 IV]$ 
and $[III, IV]$ and their quadratic twists exist if and only if $p\equiv 5\mod 6$.
\end{Corollary}

\emph{Proof:} In the first case, a quadratic twist of $S(X)$ has one of the two given additive configuration. This gives a contradiction by the previous corollary. 
The same argument applies to the second case if $p\equiv 1\mod 6$. Conversely,  
if $p\equiv 5\mod 6$, the existence of counterexamples with the given 
additive configuration is easily derived from the elliptic surfaces $Y, \widetilde Y$ and $\widehat Y$.
\qed

\section{Uniqueness of maximal counterexamples}
\label{s:minimising}

We call a counterexample for $DS(M)$ mod $p$ maximal if it attains the bounds of Thm.~\ref{Thm:p} (d). (The terminology refers to the corresponding fibre at $\infty$ being maximal.) In Rem.~\ref{Rem:sharp}, we have seen that for fixed $p$ there are infinitely many maximal counterexamples for both, odd and even $M$. These arise from rational elliptic surfaces with configuration $I_2^*, II, II$ (i.e.~$\widetilde Y$) resp.~$I_1, III^*, II$ (i.e.~$Y$) resp.~$I_6, 3 II$ via purely inseparable base change. This section aims to prove the uniqueness of any maximal counterexample with respect to $M$. Basically, this amounts to proving that the above configurations determine a unique elliptic surface up to isomorphism. In fact, we shall cover a slightly bigger class of counterexamples. 

\begin{Proposition}\label{Prop:unique}
Let $M>4$. There is a counterexample $(f,g)$ to $DS(M)$ mod $p$ with
\[
\text{deg}(f^3-g^2)=\begin{cases} 
4 \phantom{\Leftrightarrow} & 3M-2=p^r\;\;\; (\Rightarrow M \text{odd}), \\ 
5  \;\;\;\Leftrightarrow  & 6M-5=p^r,\\
6  \phantom{\Leftrightarrow} & M-1=p^r \;\;\;(\Rightarrow M \text{even}).
\end{cases}  
\]
Moreover, any such counterexample is unique up to normalisation.
\end{Proposition}

\emph{Proof:} Let $X$ denote the elliptic surface corresponding to the pair $(f,g)$. By Cor.~\ref{Cor:sep}, $X$ is not separable, so we shall also consider $S(X)$. By Lem.~\ref{Lem:s-s}, both $X$ and $S(X)$ have at least two additive fibres away from $\infty$. Since the degree of the discriminant is small by assumption, there are no multiplicative singular fibres away from $\infty$ by Remark \ref{Rem:multiplicity}. In the first step, we use this information to prove that $S(X)$ is a rational elliptic surface and to determine the configuration of singular fibres on $S(X)$. This will establish the equivalence statement of the proposition. The second step proves the uniqueness of these rational elliptic surfaces. 

\emph{1. step: Configuration on $S(X)$}

If deg$(f^3-g^2)=4$, then $M$ is odd by Thm.~\ref{Thm:p} (d). Hence $X$ has configuration $I_{e-10}^*, II, II$. Consider the quadratic twist $X'$ over $\infty$ and one $II$ fibre. Then deg $\NN_{X'}=5$, so $S(X')$ is a rational elliptic surface. It follows that $S(X')$ has singular fibres $I_2, IV^*, II$. Since purely inseparable base change commutes with quadratic twisting, $S(X)$ has configuration $I_2^*, II, II$. In particular such a surface exists and $e-10=2\,p^d$, so the claim follows from $e=6M+6$.

Let deg$(f^3-g^2)=5$. If $M$ is even, $X$ has a fibre $I_{e-5}$. We have seen in the proof of Lem.~\ref{Lem:4,5}, that this implies the configuration of $S(X)$ to be $I_1, III^*, II$ and $e-5=6M-5=p^d$. On the other hand, if $M$ is odd, then $X$ has configuration $I_{e-11}^*, III, II$. The argument from the degree $4$ case above shows that $S(X)$ has configuration $I_1^*, III, II$. In particular, $e-11=6M-5=p^d$.

Assume deg$(f^3-g^2)=6$. If $M$ is even, then $X$ has a fibre $I_{e-6}$ and furthermore
\[
IV, II\;\;-\;\; III, III\;\;-\;\; 3II.
\]
Since $e-6\equiv 6\mod 12$, $S(X)$ has a fibre of type $I_{6n}$ with $n$ odd. Using the degree of the conductor of $X$, Thm.~\ref{Thm:P-S} allows only four configurations on $S(X)$
\[
I_6, IV, II\;\;-\;\; I_6, III, III\;\;-\;\; I_6, 3II\;\;-\;\;I_{18}, 3II.
\]
All but the third configuration can be ruled out by quadratic twisting. On the other hand, an elliptic surface $\widehat{Y}$ with the third configuration implicitly appeared in Rem.~\ref{Rem:sharp}, realised as a base change of $\widetilde Y$ via $\pi_3$. In particular, $e-6=6M-6=6\,p^d$.

If $M$ is odd, then $X$ has a fibre $I_{e-12}^*$. Hence $S(X)$ has a fibre $I_{12n}^*$ with $n>0$. In particular, $S(X)$ is not rational. Reasoning as above leads to the following configurations on $S(X)$:
\[
I_{12}^*, IV, II\;\;-\;\; I_{12}^*, III, III\;\;-\;\; I_{12}^*, 3II\;\;-\;\;I_{24}^*, 3II.
\]
All configurations can be ruled out by quadratic twisting.

\emph{2. step: Uniqueness of rational elliptic surfaces}

We want to show that a counterexample $(f,g)$ as above is unique up to normalisation. To achieve this, we shall work with the corresponding separable rational elliptic surface $S(X)$. The equivalent statement then is that $S(X)$ is uniquely determined by its configuration up to isomorphism. Since $I_1^*, III, II$ is a quadratic twist of $I_1, III^*, II$, we only have to consider the following three configurations:

\begin{Lemma}\label{Lem:config}
Let $k$ be an algebraically closed field of characteristic $p\geq 0, p\neq 2,3$. Then there is a unique elliptic surface $Z$ over $\PP_k^1$ with a section and configuration $I_1, III^*, II$ resp.~$I_2^*, II, II$ resp.~$I_6, 3 II$.
\end{Lemma}

\emph{Proof:} In characteristic $\neq 2,3$, we can always argue using the corresponding $J$-map. Over $\C$, its uniqueness can be derived from Grothendieck's theory of dessins d'enfant in all three cases. In general, however, we have to work with explicit equations. Using M\"obius transformation, we can normalise three cusps. It is immediate that this gives the uniqueness for the first two configurations (thus corresponding to $Y$ and $\widetilde Y$). 

For the remaining configuration, we mentioned the elliptic surface $\widehat{Y}$. After a quadratic twist, this has Weierstrass equation
\[
\widehat{Y}:\;\;y^2 = x^3 - 3\,t\,(t^3-1)\,x + (t^3-1)(2\,t^3-1).
\]
Hence we shall assume the additive fibres of $Z$ to sit above the third roots of unity. Then the Weierstrass form of $Z$ reads
\[
Z:\;\;\; y^2 = x^3 - 3\, (t^3-1)\, A\, x + 2 \,(t^3-1)\, B
\]
with $A$ linear and $B$ of degree at most $3$. We assume that the discriminant of $Z$ is
\[
\Delta = 1728\, (t^3-1)^2\, [(t^3-1)\, A^3 - B^2] = - 1728\, (t^3-1)^2\, C^6
\]
for some linear polynomial $C$, possibly constant, non-zero at the third roots of unity. A substantial simplification is achieved by evaluating $\Delta$ at the prescribed cusps. At a third root of unity $\rho$, this gives
\[
B(\rho) = \pm C(\rho)^3.
\]
Each choice of sign results in three linear equations for the coefficients of $B$. Using these simple relations, we continue by comparing coefficients in $\Delta$. One easily derives that there are exactly four solutions. These can be identified with each other using the variable changes $t\mapsto \frac 1t$ (exchanging $0$ and $\infty$) and $(x,y) \mapsto (-x,\sqrt{-1} y)$ (changing the sign of $B$). \qed

We finish this section with a slightly different uniqueness result.
Recall that by Criterion \ref{Crit} there are no counterexamples to $DS(M)$ mod $p$
if $p>6M-5$. 

\begin{Lemma}
Every counterexample to $DS(M)$ mod $p$ with $p>3M$ arises via variable
transformation, purely inseparable base change and quadratic twisting from the surface
$Y$ in eq.~(\ref{eq:Y}).
\end{Lemma}

\emph{Proof:} Consider the corresponding elliptic surface $X$. By Rem.~\ref{Rem:multiplicity}, every fibre of type $I_n^{(*)}$ has index a multiple of $p$. Let $mp$ be the sum of these indices. By 
Lem.~\ref{Lem:s-s}, there are at least two other additive fibres. Hence $e(X)\geq mp+4$. Independent of the question whether the fibre at $\infty$ is reduced, we obtain 
\[
M\geq \frac 16 (mp+4).
\]
The assumption $p>3M$ gives $p>\frac m2 p+2$, so $m=1$. Thus $S(X)$ has exactly one fibre of type $I_n^{(*)}$ with $n>0$, and moreover $n=1$. We shall now use the fact that in characteristic $p\neq 2,3$, every elliptic surface can be realized up to quadratic twisting as pull-back from $Y$ via its $J$-map. In consequence, the sum of the indices of the $I_n^{(*)}$ fibres equals the degree of $J$. In the present case of $S(X)$, it follows that $J$ has degree $1$. Hence it is a M\"obius transformation. \qed

\section{Integral points on $\mathcal E:\;\; Y^2 = X^3 + (t^2-1)^2$}
\label{s:integral}

Following a famous idea of \v Safarevi\v c, Shioda in \cite{Sh-DS} emphasised the interplay between the elliptic surfaces $X_{f,g}$ and $X_{0,g^2-f^3}$: The pair $(f,g)$ which gives the former, is a section of the latter. Motivated by this, we shall study the following elliptic curve $\mathcal E$ over $k(t)$ which arises from $\widetilde Y$:
\[
\mathcal E:\;\;\; Y^2 = X^3 + (t^2-1)^2.
\]
Here $k$ denotes an algebraically closed field of characteristic $p\geq 0, p\neq 2,3$ as before. 
We are particularly concerned with the integral points of $\mathcal E$, i.e.~$(f,g)\in \mathcal E(k[t])$. Since consequently  $f^3-g^2=-(t^2-1)^2$, this question is closely related to the theory of Davenport-Stothers triples. 

Denote by $G$ the group of group automorphisms of $\mathcal E$ (with respect to the group law). Then 
$$
\begin{array}{ccccl}
G = <\varphi, \tau>  \;\; \text{ with } &  \varphi:  & (X,Y) & \mapsto & (\rho X, -Y)\;\;\;\;(\rho^3=1, \rho\neq 1),\\
& \tau: &  t & \mapsto & \;\;\;-t.
\end{array}
$$
The group $G$ operates on the integral points of $\mathcal E$. The following proposition gives all integral points on $\mathcal E$ in terms of representatives up to the $G$-action. 

\begin{Proposition}\label{Prop:integral}
In the above notation, let
$$
\begin{array}{lcl}
P=(t^2-1, t(t^2-1)), && Q=(2 (t+1), (t+1)(t+3)),\\
R=(0, (t^2-1)), && S=(-\frac{t^2+3}3, \frac{t(t^2-9)}{3\sqrt{-3}}).
\end{array}
$$
Then $\mathcal E(k[t]) = G \cdot \{P, Q, R, S, Frob_q^* P; \; q=p^r\equiv 1\mod 6\}$.
\end{Proposition}

\emph{Proof:} Let $(f,g)$ be an integral point of $\mathcal E$. Consider the elliptic surface $X=X_{f,g}$ corresponding to this point. This has discriminant $\Delta=1728\,(f^3-g^2)=-1728\,(t^2-1)^2$, so the fibres above $\pm 1$ have type $I_2$ or $II$.

If deg $f\leq 4$ and deg $g\leq 6$, then $X$ is rational. Hence the fibre at $\infty$ has type $I_8, I_2^*$ or $IV^*$. In total, this gives nine possible configurations. We shall now sketch how to rule out five of them.

Start with the configuration $IV^*, 2 I_2$. Apply a cubic base change which ramifies exactly at $1$ and $\infty$. The resulting configuration is $I_6, 3 I_2$. A corresponding elliptic surface would again be rational. Hence, the N\'eron-Severi lattice would have discriminant $-1$. With the above configuration, this is impossible. 

Having ruled out this configuration, we can exclude four further configurations by quadratic twisting and Thm.~\ref{Thm:P-S}. The remaining four configurations do in fact exist - they correspond to integral points of $\mathcal E$ as follows:
$$
\begin{array}{rcr}
I_2^*, 2 II  -  P &\;\;\;& IV^*, II, I_2  -  Q\\
IV^*, 2 II  -  R && I_2^*, 2 I_2  -  S
\end{array}
$$
Here we only have to take two easy facts into account: Each configuration determines a unique elliptic fibration with section up to isomorphism. After normalising the cusps and scaling, every isomorphism of fibrations is given by an element of $G$. 

On the other hand, if deg $f>1$ or deg $g>2$, then deg $f=2M$, deg $g=3M$ due to the low degree of $\Delta$. Hence the theory of Davenport-Stothers triples applies. In particular, if $M>3$, then we have seen in the proof of Prop.~\ref{Prop:unique} (Thm.~\ref{Thm:p} (d)) that $X_{f,g}$ arises from the rational elliptic surface corresponding to $P$ by base change via $Frob_q, q=p^r\equiv 1\mod 6$. In other words, $(f,g)=Frob_q^* P$ as claimed.

Finally, if $M=3$, then $X$ is a K3 surface with an $I_{14}^*$ fibre at $\infty$. By Thm.~\ref{Thm:P-S}, the other two singular fibres have type $II$. Consider the quadratic twist $X'$ at $1$ and $\infty$. By our previous argumentation (cf.~proof of Prop.~\ref{Prop:unique}), $X'$ is not separable. Thus $X$ is not separable either. It follows that $p=7$ and $(f,g)=Frob_7^* P$. \qed

\begin{Remark}
In view of Prop.~\ref{Prop:unique} and its proof, one can achieve a similar result for the elliptic curves $\mathcal E':\;\; Y^2 = X^3 + t^3(t-1)^2$ \;and\;\; $\mathcal E'':\;\; Y^2 = X^3 + (t^3-1)^2$.
\end{Remark}

\section{Characteristic 3}
\label{s:3}

In characteristic 2 and 3, the polynomial relation from Theorem
\ref{Thm:DS} is pointless (cf.~\cite[Rem.~on p.~55]{Sh-DS}). Similarly, the 
Weierstrass form (\ref{eq:Weier}) is not available in these characteristics. 
However, we can still ask whether there
are elliptic surfaces over $\PP^1$ with fibres which were
impossible over $\C$ (for fixed Euler number).  This question is
subtle due to the possibility of wild ramification. For
rational and K3 elliptic surfaces, the answer is negative
(cf.~\cite{S-max}). For semi-stable fibrations, the same follows in full generality from Rem.~\ref{Rem:s-s-2,3}.

This section concerns characteristic $3$. We shall start with an example: Consider the extremal rational elliptic surface
\begin{eqnarray}\label{eq:X}
X/\F_3:\; y^2 = x^3 + t^2 x^2 + t x + 1.
\end{eqnarray}
This has singular fibres $I_9$ and $II$ with wild ramification at
the latter and a section $(0,1)$ of order $3$. Let $X_3$ be the
pull-back of $X$ via $\phi_3$. It follows from Tate's algorithm \cite{Tate} that its singular fibre at $0$ has type $IV^*$. In
particular, $X_3$ is extremal and unirational (thus supersingular). 

On the other hand, we can transfer the * to the cusp at $\infty$
to obtain $X_3^*$ with singular fibres $I_{27}^*, II$. As
$e(X_3^*)=36$, the $I_{27}^*$ fibre contradicts the $\C$-bounds for the singular fibres. Using other cyclic base changes $\phi_d$ as well, we will derive the

\begin{Proposition}\label{Prop:3}
Let $e\in 12\,\Z$ with $e\geq 96$. In characteristic $3$, there are elliptic
surfaces over $\PP^1$ with Euler number $e$ contradicting the
$\C$-bounds for multiplicative and additive fibres (that is, with
a fibre of type $I_n, n\geq \frac 56 e$, resp.~$I_m^*, m\geq
5(\frac e6-1)$).
\end{Proposition}

The proof of the proposition mimics that of Lemma \ref{Lem:d_0} and Corollary
\ref{Cor:M_0} for the following base changes:
\[
X_d \stackrel{\phi_d}{\longrightarrow}X \;\;\; \text{ and their
*-twists } X_d^*
\]
where we transfer the * to $\infty$ (as for $X_3^*$ above) or add two $^*$s above $0$ and $\infty$. This method requires an
estimate for the Euler number $e(X_d)$. This is provided by the
following

\begin{Lemma}\label{Lem:wild}
Let $X_d$ be as above. Then the wild ramification $\delta$ at the fibre
above $0$ and the Euler number of $X_d$ are given by the following two cases:
\begin{enumerate}[(i)]
 \item Let $r\in\N$ and $d=3^r$. Then $\delta=1$. In particular, the fibre at $0$ has type $IV^*$ or $II$, according to \[e(X_d)=9d+\begin{cases}9 &
 \text{if $r$ is odd},\\ 3 & \text{if $r$ is even.}\end{cases}\]
 \item Let $(d,3)=1$. Then $\delta=d$. In particular, the fibre at $0$
 behaves as in characteristic $0$ under the base change:

$$
\begin{array}{|c||c|l|}
\hline
d\mod 6 & \text{fibre type} & \,\;\;e(X_d)\\
\hline
\hline
1 & II & 10d+2\\
2 & IV & 10d+4\\
4 & IV^* & 10d+8\\
5 & II^* & 10d+10\\
\hline
\end{array}
$$

\end{enumerate}
\end{Lemma}
We verified both statements explicitly using Tate's
algorithm. They can also be derived from \cite[Thm.~2.1]{MT}.

Combining the two statements of Lemma \ref{Lem:wild} for arbitrary
$\phi_d$, we obtain the estimate
\[
e(X_d)\leq 10(d+1).
\]
Then we apply the techniques of Section \ref{s:counterexamples} to find that
Proposition \ref{Prop:3} is guaranteed if $d\geq 18$ (that is,
$e\geq 180$). Going backwards and checking base changes in detail,
we verify that there are in fact enough examples to cover all
Euler numbers from $96$.\qed

\begin{Remark}
For $e=12, 24$, the $\C$-bounds hold by \cite{S-max}. By Prop.~\ref{Prop:3} plus some further explicit surfaces $X_d^{(*)}$, the only fibre types and Euler numbers where the $\C$-bounds also might hold, are multiplicative for $e=36,48$, and
additive for $e=72, 84$.
\end{Remark}

\begin{Corollary}
Let the characteristic be $3$.  An elliptic surface over $\PP^1$ with Euler number $e$
has a fibre of type $I_{e-3}$ or $I_{e-9}^*$, if and only if it is isomorphic to $X_{3^{2r}}$
or $X_{3^{2r-1}}^*$ for some $r\in{\N}_0$. Fibres of type $I_{e-4}$ and $I_{e-8}^*$ are impossible for $e>12$.
\end{Corollary}

\emph{Proof:} Let $Z$ be as in the Corollary. In any case, deg $\NN\leq 5$. Hence, $S(Z)$ is rational. The restrictions now follow as in the proof of Lemma~\ref{Lem:4,5}, although we have to take wild ramification into account. For instance, $e-4$ is coprime to $3$, so a surface with an $I_{e-4}$ fibre is separable. Hence it is rational. The same argument applies to the fibre of type $I_{e-8}^*$.

For the remaining two fibre types, consider $S(Z)$. After transfering the $^*$, if necessary, we obtain a surface $S'$ with two singular fibres, one multiplicative and one additive with wild ramification of index $1$. In particular, $S'$ is extremal. Extremal rational elliptic surfaces with at most three singular fibres have been classified by W.~E.~Lang in \cite{L2}. We find three surfaces as above, each uniquely determined by the singular fibres up to isomorphism. These surfaces are $X$ from (\ref{eq:X}) and surfaces with configuration $IV^*, I_3$ and $II^*, I_1$. As only the last one is Frobenius-minimal, $S'$ equals this surface by construction. In particular, $X$ is obtained from $S'$ via $\phi_9$, so the claim follows. \qed

We conclude the investigation of characteristic $3$ with the remark that the surfaces $X_{3^r}$ and $X_{3^r}^*$ are unirational by construction.

\section{Characteristic 2}
\label{s:2}

In characteristic $2$, a similar result can be obtained for
multiplicative fibres by the same means:

\begin{Proposition}\label{Prop:2}
Let $e\in 12\,\Z$ with $e>24, e\neq 60$. In characteristic $2$,
there is an elliptic surface over $\PP^1$ with Euler number $e$,
contradicting the $\C$-bounds for multiplicative fibres (that is,
with a fibre of type $I_n, n\geq \frac 56 e$).
\end{Proposition}

We omit the proof. It imitates the argumentation of Proposition
\ref{Prop:3}, starting from the extremal elliptic surface with
singular fibres $I_8, III$ (cf.~\cite{L2}). 

This approach produces a unirational elliptic surface with a fibre of type $I_{e-4}$ whenever $e-4=2^{2r+1}$. Note that type $I_{e-3}$ with $e>12$ is again impossible by the proof of Lemma~\ref{Lem:4,5}, since then the elliptic surface is separable.

In search for the maximal additive fibres in characteristic $2$, one can prove 
the following analog of Lem.~\ref{Lem:4,5}: 

\begin{Lemma}
In characteristic $2$, fibre types $I_{e-8}^*, I_{e-9}^*$ and $I_{e-10}^*$ imply the rationality of an elliptic surface over $\PP^1$; in other words, they
are only possible for $e=12$.
\end{Lemma}

The techniques involved in the proof of this lemma make extensive use of Tate's algorithm \cite{Tate}. Thus they differ substantially from the approach employed throughout the present paper. Therefore we decided to postpone this issue for a future work. Beyond this small result, hardly anything seems to be known. Turning to quasi-elliptic fibrations, it is always possible to give a fibration over $\PP^1$ with a fibre of type $I_{e-8}^*, e\in 12\N$:
\begin{eqnarray}\label{eq:W_e}
W_e:\;\; y^2 = x^3 + t^2 x + t^{\frac e2 -1}.
\end{eqnarray}

\section{N\'eron-Severi group over $\F_p$}
\label{s:NS}

This section concerns the field of
definition for generators of the N\'eron-Severi group. Consider
the rational elliptic surface $X$ over $\F_3$ introduced
in Sect.~\ref{s:3}. It is immediate that $NS(X)$ is generated over
$\F_3$. The same holds for the pull-backs
$X_{3^{r}}$. It can also be checked for $X_{3^{2r-1}}^*$.  We find this worth noticing since
it is contrary to the situation for supersingular K3 surfaces (for odd characteristic, see \cite[(6.8)]{A}; for $p=2$, cf.~\cite[Thm.~4.4]{S-arith}). To understand this, we need the following result which resembles the Artin invariant for supersingular K3 surfaces:

\begin{Theorem}
\label{Thm:NS}
Let $X$ be a smooth projective surface over a finite field $k$ of
characteristic $p$. Assume that $X$ is supersingular in Shioda's sense, i.e.~$\rho(X)=b_2(X)$. Let
$N=NS(X)/NS(X)_\text{tor}$. Then
\[
\text{discr}\, N =
(-1)^{b_2(X)+1}\,p^{2\sigma}\;\;\;\;(\sigma\in\N[0]).
\]
\end{Theorem}
C.~Liedtke pointed out that Thm.~\ref{Thm:NS} can be derived from results of Illusie \cite[5.21]{I}.
Here, we sketch an alternative proof: 

For any $\ell\neq p$, we have an isomorphism
\[
N\otimes_{\Z}{\Z}_\ell \cong H^2(X, \Z_\ell).
\]
This implies that the discriminant equals $\pm p^r$. The above sign follows from the Algebraic Index Theorem. The evenness of the exponent $r$ is seen as follows:

Let $k'=\F_q$ be a finite extension of $k$ such that $NS(X)$ is generated by divisors over $k'$. Without loss of generality, assume $\F_{p^2}\subset k'$. By construction, the Tate Conjecture \cite{Tate-C} holds for $X/k'$. This is equivalent to the Artin-Tate Conjecture \cite[(C)]{Artin-Tate} by \cite{Milne} (for characteristic $2$, cf.~\cite[p.~25]{Milne-A}). Hence
\begin{eqnarray}\label{eq:A-C}
q^{\alpha(X)} |NS_\text{tor}(X)|^2 = |Br(X/k')|\; |\text{discr } N|
\end{eqnarray}
where $\alpha(X)=\chi(X) -1\, +$ dim Pic Var$(X)$.  By assumption, $q$ is a square. The same holds for $|Br(X/k')|$ by \cite{Brauer}. Hence the claim follows. \qed

With this theorem, we can directly generalise Artin's argument for (supersingular) K3 surfaces (where Artin used $\alpha=1$):

\begin{Lemma}\label{Lem:Art}
Let $k$ be a finite field and $X/k$ be a smooth projective surface. If $\F_{p^2}\not\subset k$ and $\alpha(X)$ is odd, then 
\[
\rho(X/k) < b_2(X).
\]
\end{Lemma}

\emph{Proof:} Assume that $\rho(X/k)=b_2(X)$. Then the Artin-Tate Conjecture applies to $X/k$. We employ the notation of the proof of Thm.~\ref{Thm:NS} (with $q=\# k$). In eq.~(\ref{eq:A-C}), the right hand-side is a square by \cite{Brauer} and Thm.~\ref{Thm:NS}. Hence, this also holds for $q^{\alpha(X)}$. This is equivalent to the claim. \qed

Further surfaces with $\rho(X/\F_p)=b_2(X)$ can easily be constructed:
Over $\F_2$, pull-back from the rational elliptic surface
with fibres $I_8, III$ via $\phi_{2^{2r}}$ as in the previous section. Moreover, Tate's algorithm \cite{Tate} shows that the exceptional fibre of the quasi-elliptic surface $W_e$ from (\ref{eq:W_e}) has all components of the singular fibre defined over $\F_2$ if and only if $e\in 12\N\setminus 24\N$.

Secondly, let $p>3$ and consider the surface $Y$ from (\ref{eq:Y}) over $\F_p$. Let $Y'$ denote the quadratic twist of $Y$
over $\F_p(\sqrt{3})$. Then the fibre of $Y'$ at $\infty$ has split multiplicative reduction. Thus, all fibre components of the pull-back via $\phi_{p^{r}}$ are defined over $\F_p$. Hence, the pull-back has N\'eron-Severi group generated over
$\F_p$ if and only if the non-trivial section of $Y'$ is defined over $\F_p$. Since for $Y$, this section is given by $(\frac 23 t^2, \sqrt{\frac 8{27}}t^3)$, and we twist over $\F_p(\sqrt{3})$, the last condition is equivalent to
$\sqrt{2}\in\F_p$. 

\begin{Lemma}
Let $p=2, 3$ or $p\equiv\pm 1\mod 8$. There are surfaces $X/\F_p$
(elliptic over $\PP^1$) with unbounded Euler number such that
\[\rho(X/\F_p)=b_2(X).\]
\end{Lemma}

\vspace{0.8cm}

\textbf{Acknowledgement:} 
This paper benefitted greatly from discussions with B.~van Geemen and T.~Shioda while the first author enjoyed the hospitality of Dipartimento di Matematica "Frederico Enriques" of Milano University. Support from the CRTN-network "Arithmetic Algebraic Geometry" is gratefully acknowledged. During the revision and extension of the paper, the first author was generously funded by DFG under grant "SCHU 2266/2-1". We also thank I.~Bouw for explanations about reduction properties of base changes, and C.~Liedtke for pointing out the reference to Illusie's paper.

\vspace{0.8cm}

\begin{center}
\begin{tabular}{llll}
Matthias Sch\"utt & \hspace{1cm} && Andreas Schweizer\\
Mathematics Department &&& Mathematical Sciences\\
Harvard University &&& University of Exeter\\
Science Center &&& Harrison Building\\
1 Oxford Street &&& North Park Road\\
Cambridge, MA 02138 &&& Exeter EX4 4QF\\
USA &&& United Kingdom\\
{\tt mschuett@math.harvard.edu} &&& {\tt A.Schweizer@exeter.ac.uk}
\end{tabular}
\end{center}

\end{document}